\documentclass{amsproc}
\def\card{{\mathrm{card}\,}}
\def\area{{\mathrm{area}\,}}
\def\Z{{\mathbf{Z}}}
\def\R{{\mathbf{R}}}
\def\bR{{\mathbf{\overline{R}}}}

\theoremstyle{definition}

\theoremstyle{remark}

\numberwithin{equation}{section}

\begin{document}

\def\R{{\mathbf{R}}}
\def\Sp{{\mathbf{\overline{C}}}}
\def\C{{\mathbf{C}}}
\def\U{{\mathbf{U}}}
\def\diam{{\mathrm{diam}\,}}
\def\id{{\mathrm{id}}}
\title{\bf Geometric theory of meromorphic functions}
\author{A. Eremenko}
\date{\today}
\address{Department of Mathematics, Purdue University, West Lafayette,
Indiana 47907}
\email{eremenko@math.purdue.edu}
\thanks{Supported by NSF grant DMS 0555279 and 0244547.}
\dedicatory{Dedicated to A. A. Goldberg.}
\begin{abstract}
This is a survey of results on the following problem. Let $X$ be a simply
connected Riemann surface spread over the Riemann sphere. 
How are the 
properties of the uniformizing function of this surface related to
the geometric properties of the surface?   
Based on the lecture in U. Michigan in May 2006. 
\end{abstract}

\maketitle

{\bf 1}. According to the Uniformization Theorem, for every
simply connected Riemann
surface $X$ there exists a conformal homeomorphism $\phi:X_0\to X$, where
$X_0$ is one 
of the three standard regions,
the Riemann sphere $\Sp$, 
the complex plane $\C$ or 
the unit disc $\U$. 
We say that the {\em conformal type}
of $X$ is {\em elliptic,
parabolic}
or {\em hyperbolic}, respectively. 
The map $\phi$ is called the {\em uniformizing map}.
If $X$ is given by some geometric construction, the problem arises
to relate
properties of $\phi$ to those of $X$. This includes the determination of
the conformal type of $X$ \cite{Doyle,Nev1,Volk,Wit}.
The case which was studied most is
that $X\subset\C$ is a simply
connected region, $X\neq\C$. Then $X$ is of
hyperbolic type and $\phi$ is a univalent function in $\U$. An example
of the result relating geometric
properties of $X\subset\C$ and properties of $\phi$
is the classical theorem of Caratheodory: $\phi$ is continuous in
the closed disc if and only if $\partial X$ is locally connected.
Another example is the Ahlfors distortion theorem which relates
the growth of a uniformizing map to the geometry of the image domain.

We recall a more general construction of $X$.
A surface {\em spread over the sphere}
is a pair $(X,p)$, where $X$ is a topological surface and $p:X\to\Sp$
a continuous, open and discrete map.
This map $p$ is usually called the {\em projection}.
The natural equivalence relation is $(X,p)\sim(Y,q)$ if
there is a homeomorphism $\phi:X\to Y$ with the property $p=q\circ \phi$.
According to a theorem of Sto\"{\i}lov,
every continuous open and discrete map $p$ between surfaces
locally looks like $z\mapsto z^n$. Those points where $n>1$ are
isolated, they are called {\em critical points}, or multiple points
of multiplicity $n$.
Sto\"{\i}lov's
theorem implies that
there is a unique conformal structure on $X$ which makes $p$ holomorphic.
If $\phi$ is a
uniformizing map, then $f=p\circ\phi$ is a meromorphic function in
one of the three standard regions $\Sp,\C$ or $\U$. The surface $(X,p)$
spread over the sphere
is then the ``Riemann surface of $f^{-1}$''. We also call it the
surface associated to $f$, or say that $f$ is associated with $(X,p)$. 

If $D$ is a region on the sphere, a {\em branch} of $p^{-1}$ in $D$
is a continuous function $\psi:D\to X$ such that $p\circ\psi=\id_D$.

We can define the length of a curve in $X$
as the spherical length\footnote{We choose the spherical length element to be
$2|dz|/(1+|z|^2)$, so that the curvature of the spherical metric is $+1$.}
of its image under $p$. Then $X$ becomes a metric space with
an {\em intrinsic metric}, which means that the distance between two points
is the infimum of the lengths of curves connecting these points.
Similarly, if $p:X\to\C$,
and the Euclidean metric in $\C$ is used to measure
lengths of curves,
we obtain a Riemann surface {\em spread over the plane}.

The intrinsic metric on $X$ is a smooth Riemannian metric of constant
curvature on the
complement of the critical set of $p$. 
It is easy to show that the intrinsic metric on $X$ determines
the projection $p$ up to an isometry of the image sphere or the plane.
In what follows, unless otherwise stated,
$(X,p)$ denotes a simply connected surface spread over
the sphere, equipped with the intrinsic spherical metric.

Some criteria of conformal type can be stated in terms of topological
(or even set-theoretic)
properties of $p$. 
For example, Picard's theorem implies that $X$ is of hyperbolic type
if $p$ omits three points.

The following result is due to Nevanlinna.
A point $a\in\Sp$ is called a {\em totally ramified value}
(of multiplicity $m\geq 2$)
of a surface $(X,p)$ spread over the sphere, if all preimages of $a$
under $p$ are multiple, (of multiplicity at least $m$). We allow
$m=\infty$ which means that the value $a$ is omitted.
\vspace{.1in}

\noindent
1.1 {\sc Nevanlinna's theorem} {\em If a surface spread over the sphere
has $q$ totally ramified values of multiplicities $m_k,\, 1\leq k\leq q$,
and
$$\sum_{k=1}^q\left(1-\frac{1}{m_k}\right)>2,$$
then the surface is hyperbolic. In particular, a parabolic surface
has at most four totally ramified values.}
\vspace{.1in}

The example $(\C,\wp)$ shows that a parabolic surface can indeed have $4$
totally ramified values.
Moreover, one can easily show that the equation
$$\sum_{k=1}^q\left(1-\frac{1}{m_k}\right)=2$$
has $5$ solutions up to a permutation of the $m_k$,
and to each of these solutions corresponds an elliptic or trigonometric
function. A very short proof of Theorem 1.1 can be found
in \cite{Robin}.
\vspace{.1in}

\noindent
1.2 {\sc Ahlfors's Five Islands Theorem.} {\em Suppose that for
five
Jordan regions with disjoint closures on the sphere,
there are no branches of $p^{-1}$ in any of these regions.
Then $X$ is of hyperbolic type.}
\vspace{.1in}

This theorem was stated for the first time
by Bloch \cite{Bloch1} (with discs instead of Jordan regions)
and proved by Ahlfors in
\cite{Ahl1}, as a corollary from his ``Uberlagerungsflachentheorie''.
It is a recent discovery \cite{Berg} that actually 
Theorem 1.2 can be derived from Theorem 1.1 by a simple
argument.
This derivation is based on the following important principle:
\vspace{.1in}

\noindent
1.3 {\sc Zalcman's Lemma}. {\em Let $F$ be a family of meromorphic
functions in some region $D$, which is not normal in $D$.
Then there exists a sequence $f_n\in F$, and two sequences
$r_n>0$ and $z_n\in D$ such that there exists a non-constant limit
$$f(z)=\lim_{n\to\infty}f_n(r_nz+z_n),$$
uniform on compact subsets of $\C$, and moreover,}
$$f^{\#}(z):=\frac{|f'(z)|}{1+|f(z)|^2}\leq f^\#(0)=1,\quad z\in\C.$$

{\em Proof.} Without loss of generality we may assume that $D$ is the
unit disc, functions of the family $F$ are meromorphic in the closure
of $D$, and $F$ is not normal at $0$.
This means that there exists sequence $f_n\in F$ and $w_n\to 0$ such that
$f^\#_n(w_n)\to\infty$. Then
$$\max_{z\in D}(1-|z|)f_n^\#(z)=(1-|z_n|)f_n^\#(z_n)
\geq(1-|w_n|)f^\#_n(w_n)\to\infty.$$
Thus $r_n=1/f^\#_n(z_n)=o(1-|z_n|)$. We claim that $r_n$
and $z_n$ have the required property.
Indeed, putting $g_n(z)=f_n(r_nz+z_n)$, we obtain
$g_n^\#(0)=1$ and
$$g_n^\#(z)=r_nf^\#_n(r_nz+z_n)\leq r_n
f^\#_n(z_n)\frac{1-|z_n|}{1-|r_nz+z_n|}\leq\frac{1-|z_n|}{1-|z_n|-r_n|z|}
\to 1,$$
because $r_n=o(1-|z_n|)$.
So $g_n$ is a normal family. After selecting a subsequence we get
$g_n\to f$ for some meromorphic function $f$,
$f^\#(z)\leq 1,z\in\C$ and $f^\#(0)=1$,
which proves all statements of the lemma.
\vspace{.1in}

Now we derive Theorem 1.2 by contradiction. Suppose that
$f$ is a meromorphic function in the plane, and $D_k$ are five
disjoint Jordan regions such that there are no inverse branches
of $f^{-1}$ in $D_k,\, 1\leq k\leq 5$.
Choose five points $a_k$ on the sphere and consider quasiconformal
homeomorphisms $\psi_n$ of the sphere which map each $D_k$ into
$1/n$-neighborhood of $a_k$. Then the surfaces $(\C,\psi_n\circ f)$
spread over the sphere are all parabolic because the maps
$\psi_n\circ f$ are
{\em quasiregular}\footnote{Quasiregular maps in dimension $2$ are compositions of holomorphic functions
with quasiconformal maps. Quasiconformal maps preserve the
conformal type of a simply connected Riemann surface.}, 
so there exist homeomorphisms
$\phi_n$ of $\C$ such that $f_n=\psi_n\circ f\circ\phi_n$ are meromorphic
functions. These functions have no inverse branches over $1/n$
neighborhoods of $a_k$. We can use arbitrarity in the choice
of $\phi_n$ to normalize our functions: $f_n(0)=0,\, f'_n(0)=1$.
If the family $\{ f_n\}$ is normal in the whole plane, then the limit functions
are non-constant because of the normalization, and it is easy to see
that $a_k$ are totally ramified values of these functions, contradicting
Theorem 1.1. If $\{ f_n\}$ is not a normal family in the plane, we
apply Zalcman's lemma to make it normal, and again obtain a contradiction.
\vspace{.1in}

Zalcman's lemma shows the importance of study of meromorphic
functions with bounded spherical derivative. Very little is known
about this class. We mention a theorem of Clunie and Hayman that if
an entire function has bounded spherical derivative than it is
at most of order one, normal type. (A typical meromorphic function
of this class has order two, normal type).

Sullivan asked the general question, for which surfaces $(X,p)$
the conformal type is determined
by topological properties of $p$. More precisely, let us say that
a simply connected surface $(X,p)$ spread over the sphere has a
{\em stable type} if for every homeomorphism $\psi:\Sp\to\Sp$ the
surface $(X,\psi\circ p)$ has the same type as $(X,p)$.
For example, surfaces satisfying the conditions of Theorem 1.1 or 1.2
are of stable hyperbolic type.
Another interesting class of surfaces of stable type is the {\em Speiser
class} $S$. We say that $(X,p)\in S$ if there exists a finite set $A\subset\Sp$
such that the restriction
$p:X\backslash p^{-1}(A)\to\Sp\backslash A$
is a covering. The stability of type of such surfaces was proved
by Teichm\"uller in \cite{Teich} as one of the first applications
of quasiconformal mappings.
The argument is based on the following simple
\vspace{.1in}

\noindent
{\sc Lemma.} {\em Let $f$ be a function of Speiser class, in
$\{ z:|z|<R\},$ where $R=1$ or $R=\infty$, and $A$ is
the finite set as above. Let $f_1=\psi\circ f\circ\phi,$
where $f_1$ is a meromorphic function and $\psi$ and $\phi$
are homeomorphisms, $\phi$ leaves $0$ and $1/2$ fixed.
Suppose that there is an isotopy $\psi_t, 0\leq t\leq 1,$
$\psi_0=\id,$ $\psi_1=\psi$ such that
the elements of $A\cup\{ f(0),f(1/2)\}$
remain fixed by each $\psi_t$.
Then $f_1=f$.}
\vspace{.1in}

By the covering homotopy theorem there exists
an isotopy $\phi_t,\; \phi_1=\phi$ 
such that $\psi_t\circ f\circ \phi_t=f_1,\;
0\leq t\leq 1$. The functions $\phi_t(0)$ and $\phi_t(1/2)$
are continuous
and can take discrete sets of values, hence $\phi_t(t)(0)=0$
and $\phi_t(1/2)=1/2$ for all $t$.
Putting $t=0$ we obtain $f\circ \phi_0=f_1$.
It is easy to conclude from the last equality that
$\phi_0$ is a conformal homeomorphism of
$\{ z:|z|<R\}$. On the other hand, it fixes
$0$ and $1/2$, so $\phi_0=\id$ and $f_1=f$.
\vspace{.1in}

Now Teichm\"uller's result follows because every homeomorphism
is isotopic to a quasiconformal homeomorphism modulo any given finite set,
and quasiconformal homeomorphisms preserve conformal type of a simply
connected surface.

Teichm\"uller's argument extends to the somewhat wider class
consisting of surfaces with the property that the distances between 
their singularities\footnote{
A formal definition of singularities and distance between them
is in section 3.}
are bounded from below by a positive constant.
All surfaces of stable parabolic type known to the author have this
property.
The simplest example of a parabolic surface with a non-isolated
singularity is
$(\C,p)$, with $p(z)=\sin z/z$, and it follows from a result of
Volkovyskii \cite[Th.\ 45]{Volk} that the conformal type of this surface
is not stable.


\vspace{.1in}

{\bf 2.}
If we take the five regions in Ahlfors's theorem
to be spherical discs of equal radii,
Theorem 1.1 implies the following \cite{Ahl2}:
{\em Suppose that for some $\epsilon>0$ there are no branches
of $p^{-1}$ in discs
of radii $\pi/4-\epsilon$. Then $X$ is of hyperbolic type.}
The question arises, what is the best constant for which this
result still holds. Let $B(X,p)$ be the supremum of radii of discs where
branches of $p^{-1}$ exist, and $B=\inf B(X,p)$, where the infimum is
taken over all
surfaces of elliptic or parabolic type.
Ahlfors's estimate $B\geq\pi/4$ was improved
by Pommerenke \cite{Po} to $B\geq\pi/3$, and recently the sharp result
was obtained in \cite{BE}:
$$B=b_0:=\arccos(1/3)\approx 0.39\pi.$$
We have $B(\C,\wp)=B$, where $\wp$
is the Weierstrass function
of a hexagonal lattice. It is interesting to notice that $B=b_0$
implies Theorem 1.1 by a simple argument given in \cite{BE}.

For surfaces $(X,p)$ of elliptic type
we have $B(X,p)>b_0$, but
it is not known whether the constant $b_0$ is best possible
in this inequality.

We sketch the proof for elliptic surfaces.
First one constructs a triangulation $T$ of $X$ into geodesic triangles,
so that the vertices of this triangulation coincide with
the set of critical points, and the circumscribed
radius of each triangle is at most $B(X,p)$. This is always possible to do if
$B<\pi/2$ which we can assume. Suppose now that $B(X,p)\leq b_0$. Then the
circumscribed radius of each triangle is at most $b_0$, and an elementary
geometric argument shows that the area of each triangle is
at most $\pi$. Notice that by Gauss formula,
$\area(\Delta)=\sum\alpha(\Delta)-\pi$,
where $\alpha(\Delta)$ is the sum of the angles of $\Delta$.
As $\area(\Delta)\leq\pi$ we conclude that $\area(\Delta)\leq \alpha(\Delta)/2$.
If we denote by $\alpha(v)$ the total angle at a vertex,
then $\alpha(v)=4\pi,$
assuming that all critical points have multiplicity $2$.
If $d$ is the degree of our rational function then the
total area is
$$4\pi n=\sum_{\Delta\in T}\area(\Delta)\leq\frac{1}{2}\sum_v\alpha(v)=
2\pi(2n-2),$$
and this is a contradiction. 
\vspace{.1in}

The proof of $B\geq b_0$ 
for parabolic surfaces
is more complicated. For our class of surfaces with intrinsic metric,
one can define 
integral curvature \cite{Resh} as a signed
Borel measure on $X$ which is equal to the area on the
smooth part of $X$ and has negative
atoms at the critical points of~$p$.
The assumption that $B(X,p)<b_0$ implies that the atoms of negative
curvature are sufficiently dense on the surface, so that on large pieces
of $X$
the negative part of the curvature dominates the positive part.
Then a bi-Lipschitz modification of the surface is made, which spreads
the integral curvature more evenly on the surface, resulting in
a surface whose Gaussian curvature is bounded from above by a negative
constant, and the Ahlfors--Schwarz lemma implies
hyperbolicity. A non-technical exposition of the ideas of this proof is
given in the survey
\cite{Bonk} which contains some
further geometric applications of this technique of spreading the
curvature by bi-Lipschitz modifications of a surface.
\vspace{.1in}

{\bf 3}. To formalize the notion of a singular point of a multi-valued
analytic function, Mazurkiewicz \cite{Mazur}
introduced another metric,
$\rho(x,y)=\inf\{\diam p(C)\}$ on $X$,
where $\diam$ is the diameter with respect to the spherical metric and
the infimum is taken over all curves
$C\subset X$ connecting $x$ and $y$.
Every point $x\in X$ has a neighborhood where the Mazurkiewicz metric
coincides with the intrinsic one, but in general the Mazurkiewicz metric
is smaller. For example, on the surface $(\C,\cos)$ spread over the sphere,
the intrinsic
distance between $0$ and $2\pi k$ is $2\pi k$, while the
Mazurkiewicz distance is $\pi$.

Let $X^*$ be the completion of $X$
with respect to the Mazurkiewicz metric.
Then $p$ has a unique continuous extension to $X^*$. The elements
of the set $Z=X^*\backslash X$ are called {\em transcendental
singularities} 
of $(X,p)$. The simplest example of a transcendental
singularity is a logarithmic branch point.
The {\em algebraic singularities} are just the critical
points of $p$. The images of singularities under $p$ will
be called {\em singular values}. The images of critical points
are called {\em critical values}.

To each transcendental singularity corresponds an {\em asymptotic curve}
$\gamma:[0,1)\to X$ which has no limit in $X$ but its image
$p\circ\gamma$ has a limit in $\Sp$. This limit is called
an {\em asymptotic value} and it is the projection of the singularity. 

If $D\subset\Sp$ is a region containing no singular values then
the restriction $f:f^{-1}(D)\to D$ is a covering map.
The closure of the set of singular values is characterized by this
property. However the set of singular values need not be closed.

If $X$ is of parabolic type, then the set of singularities is totally
disconnected. This can be proved by using Iversen's theorem \cite{Iver,Nev1}.
The following classical result \cite{Nev1}, which implies Iversen's theorem,
shows that if we ``look in all directions
from a point'' on a parabolic surface then very few singularities are
visible.
\vspace{.1in}

\noindent
3.1 {\sc Gross's Theorem.}
{\em Let $(X,p)$ be a simply connected surface of
parabolic type spread over
the sphere, and $x\in X$. Then
$X$ contains a geodesic ray from $x$ of length $\pi$ (that is from $x$ to the
``antipodal point'')
in almost every direction.}
\vspace{.1in}

It is not known whether the estimate of the size of the exceptional set
in this theorem can be improved, but there are examples where this
exceptional set of directions has
the power of the continuum \cite[Th. 17]{Volk}.
We include a sketch of such example. Let $D$ be the unit disc,
and $E\subset \partial D$ a Cantor set whose complement
consists of the arcs $L_k=(a_k,b_k)$. We assume that $\partial D$
has the standard orientation which induces orientation
from $a_k$
to $b_k$ on these arcs. 
Let us cut the surface $(\C,(a_k-b_ke^z)/(1-e^z))$ along a simple arc
that projects to to $L_k$ and call $X_k$ the
piece which lies on the right of the arc. Such surfaces $X_k$ are
called {\em logarithmic ends}. Now, for each $k$ we paste
 $X_k$ to $D$
along the arc $L_k$, respecting the projection.
The result is an open simply connected surface $X$ spread over
the sphere. It is easy to show, that this surface can be parabolic
if the complementary arcs to $\cup_{k=1}^n L_k$ decrease sufficiently fast
as $n\to\infty$. To do this one approximates $X$ by surfaces with
finitely many singularities which are all parabolic by
Theorem 4.2 below. Volkovyski obtained a quantitative sufficient condition
for $X$ to be of parabolic type, but it is much stronger than the only known
necessary condition that $E$ is of zero logarithmic capacity.
Evidently, on $X$, the length of a geodesic ray
from the center of $D$ in the direction of any point of $E$
is $\pi/2$.

The projection of the set $Z$ of singular points is an analytic (Suslin) set
\cite{Mazur}, and for every analytic subset $A$ of the sphere one can find
surfaces of both parabolic and hyperbolic types for which $A=p(Z)$
\cite{Heins}.

The following classification of transcendental singularities
was introduced by Iversen
\cite{Iver}. A singular point $x\in X^*\backslash X$
is called {\em direct} if for
some neighborhood $V\subset X^*$ of $x$, the map $p$ omits $p(x)$
in $V\backslash \{ x\}$. Otherwise $x$ is called {\em indirect}.
For example, $(\C,\sin z/z)$ has two direct
singularities over $\infty$ and
two indirect singularities over $0$.
The following result was proved
in \cite{Heins1}:
\vspace{.1in}

\noindent
3.2 {\sc Heins's Theorem.} {\em For a parabolic Riemann
surface spread
over the sphere, the set of projections of direct singularities
is at most countable.}
\vspace{.1in}

On the other hand, for some parabolic surfaces, the set of direct
singularities 
lying over one
point may have the power of the continuum. 
To construct such example, one can take an infinite dyadic tree $T$
properly embedded in the plane, and consider a simply connected
neighborhood $U$ of $T$, so small that all infinite branches of $T$
define different accessible boundary points of $U$ at infinity.
Then one constructs an entire function $f$ that tends to infinity
along each infinite branch of $T$ but remains bounded in the complement
of $U$. Then the surface $(\C,f)$ has uncountably many singularities
over $\infty$, all of them direct, because $\infty$
is omitted. An explicit construction of such entire
function $f$ is given in \cite{BE7}. 


To state further results on direct singularities
we recall the notion of
the order of a meromorphic function in the plane. Let $(X,p)$
be an open simply connected
Riemann surface of parabolic type, spread over the sphere.
Then the intrinsic metric defines
a notion of area on $X$. If $\phi$ is a uniformizing map and $f=p\circ\phi$
the corresponding meromorphic function in $\C$, then the ``average
covering number'' of the sphere by the images of the discs
$D(r)=\{ z:|z|\leq r\}$ is defined as
the area of $\phi(D(r))$ divided by the area of the sphere $\Sp$, which is the same
as
$$A(r,f)=\frac{1}{\pi}\int_{D(r)}\frac{|f'|^2}{(1+|f|^2)^2},\quad r\geq 0,$$
and the order of $f$ is defined as
$$\lambda(f)=\limsup_{r\to\infty}\frac{\log A(r,f)}{\log r}.$$
It is easy to verify that the order depends only on $(X,p)$ rather then on the
choice of the uniformizing map $\phi$.
\vspace{.1in}

\noindent
3.3 {\sc Denjoy--Carleman--Ahlfors Theorem}. {\em
If $f$ is a meromorphic
function in the plane,
and the Riemann surface $(\C,f)$ has $k\geq 2$
direct singularities, then}
$$\liminf_{r\to\infty}r^{-k/2}{A(r,f)}>0.$$

If $p$ omits a point $a\in\Sp$, then a theorem of Lindel\"of implies that
at least half of all singularities of $(X,p)$ lie over $a$,
and these singularities are evidently direct. So we obtain
that such surface can have at most $2\lambda(f)$ singularities
over the points
in $\Sp\backslash \{ a\}$. As a corollary, an entire function of order
$\lambda$ can have at most $2\lambda$ finite asymptotic values.
This is in contrast with the case of meromorphic functions: there exist
meromorphic functions in the plane of arbitrary prescribed order
$\lambda\geq 0$
for which every point on the sphere is an asymptotic value \cite{Erem1}.

In the simplest example $\sin z/z$ mentioned above,
the indirect singularities over $0$ are accumulation points of
critical points. A result of Volkovyski \cite[Th.\ 17]{Volk} shows
that this is not always so: there are parabolic surfaces without critical
points, having indirect singularities. However, for functions of finite
order, the following theorem was proved in \cite{BergE}:
\vspace{.1in}

\noindent
3.4 {\sc Theorem}.
{\em Let $f$ be a meromorphic function of finite order,
and let $a$ be an indirect singularity of $(\C,f)$. Then
every neighborhood of $a$ contains critical points $z$ such that}
$f(z)\neq f(a)$.
\vspace{.1in}

This theorem can be used to prove the existence of critical points
under certain circumstances, and more generally, to study the value distribution
of derivatives. 
In \cite{BergE} it was used
to prove a result which was conjectured by Hayman:
for every transcendental meromorphic function $f$ in $\C$, the equation
$ff'=c$ has infinitely many solutions
for every $c\in\C\backslash\{0\}$. There is no growth restriction on $f$
in this last result.

The scheme of the proof of Hayman's conjecture is the following.
Let $g(z)=z-f^2/(2c)$, and suppose first that $f$ is of finite order.
Then each zero of $f$ is a fixed point
of $g$ with multiplier $1$. According to Fatou's theorem from holomorphic
dynamics, to each such fixed point a region of immediate attraction
is attached, and this region contains a singular value of $g$. By Theorem 3.4,
$g$ has to have infinitely many critical values. And each critical
point of $g$ is a solution of $ff'(z)=c$.
The case that $f$ is of finite order and has finitely many zeros is
treated separately by elementary considerations.

To extend the result to functions of infinite order one uses the
following generalization of Zalcman's lemma which is due to
X. Pang: {\em Let $F$ be a non-normal family of meromorphic
functions and $-1<k<1$. Then there exist sequences $f_n\in F,\,
r_n>0$ and $z_n$ such that the limit
$$\lim r_n^{-k}f_n(r_nz+z_n)=g$$
exists on every compact in the plane, and $g$ is a non-constant meromorphic
function with bounded spherical derivative (so the order of $g$ is
at most $2$).}

To complete the proof of Hayman's conjecture one uses $k=1/2$.
\vspace{.1in}


Goldberg and Heins independently noticed that in Theorem 3.3
one may count some indirect singularities, so-called $K$-singularities,
together with the direct ones.
However a geometric characterization of $K$-singularities
is known only for a very special class of symmetric surfaces \cite{GoldbergCV}.
\vspace{.1in}

{\bf 4.} There are few instances when a precise correspondence
can be established
between classes of surfaces spread over the sphere and classes of
meromorphic functions corresponding to them.
We mention first a class of hyperbolic surfaces, 
spread over the plane.
\vspace{.1in}

\noindent
4.1 {\sc Theorem}. {\em For a 
surface $(X,p)$ spread over the plane,
the following conditions are equivalent:
\newline
(a) The (Euclidean) radii of discs where branches of $p^{-1}$ exist are
bounded.
\newline
(b) A linear isoperimetric inequality holds on $X$.
\newline
(c) $X$ is of hyperbolic type and the uniformizing map $\phi:\U\to X$
is uniformly continuous with respect to the hyperbolic metric on $\U$
and the intrinsic (Euclidean) metric on} $X$.
\vspace{.1in}

The equivalence between (a) and (c) is essentially Bloch's theorem
\cite{Bloch}.
For the equivalence of (b) to the other two conditions the reference is
\cite{BE2} where Theorem~4.1 is stated for a more general classes of
surfaces with intrinsic metric (not necessarily spread over the plane).
Holomorphic functions
$f=p\circ\phi$, where $\phi$ satisfies (c), are called Bloch functions.
This class is important because of its connection with
univalent functions: if $g$ is univalent in $\U$, then $\log g'$ is a Bloch
function, and every Bloch function has the form $c\log g'$ where $g$ is
univalent in $\U$ and $c$ is a constant \cite{Po}.
In the case of meromorphic functions
and surfaces spread over the sphere, Ahlfors and Dufresnoy \cite{Duf} proved that
(b) implies (c). It is plausible that (b) is actually equivalent to (c)\footnote{
This was proved by Bruce Kleiner (unpublished).},
but it is not clear what sort of conditions could replace (a)
in the spherical case.
Functions satisfying (c) with the spherical metric
instead of the Euclidean one are called {\em normal}, and they were
much studied, see, for example
\cite{Pom3}.

Passing to parabolic surfaces, we first
describe a way to visualize certain classes of
surfaces spread over the sphere, in particular, the Speiser
class defined in section 1.
Suppose that all singular values belong to some Jordan curve
$\gamma$ on the sphere. We call such surfaces {\em cellular
surfaces} and the corresponding meromorphic functions {\em cellular
functions}.

The full preimage $p^{-1}(\gamma)$ defines a cell decomposition
of $X$. We recall that a cell decomposition of a surface $X$
is a representation of $X$ as a locally finite
union of disjoint sets called
{\em cells}. Each cell is a homeomorphic image
of a point (vertex), or of an open interval (edge) or of an
open disc (face). The closure of each cell is a union of this
cell with cells of smaller dimension. Notice that the boundary of
a face can contain infinitely many edges or infinitely
many vertices. The $1$-{\em skeleton}
of a cell decomposition is the union of edges and vertices.

Now consider a cellular surface $(X,p)$ spread over the sphere,
so that all singular
values belong to a Jordan curve $\gamma$.
One can show that the full preimage of $\gamma$ is a $1$-skeleton of a cell
decomposition of $X$. Its vertices of degree $2m$ are the multiple points
of multiplicity $m$ of $p$.
One can add vertices of degree $2$ by breaking an edge into two
edges, which is sometimes convenient. 
The curve $\gamma$ divides
the sphere into two regions, $D_0$ and $D_2$, which contain no
singular values. It follows that each face is mapped
on $D_1$ or $D_2$ homeomorphically. We can assign colors to the faces:
those mapped on $D_0$ are white and those mapped on $D_1$ are black.
Then our cell decomposition has the property that whenever two
faces have a common boundary edge, they are of different color.
This is consistent with the property that all vertices have
even degrees.

If $(X,p)$ is of the Speiser class, we have a finite set $A$ that contains
all singular values. 
Let $\gamma$ be a Jordan curve in $\C$ passing through all points of
the set $A$. Consider the cell decomposition of the sphere
whose $2$-cells are the two complementary domains to $\gamma$,
$0$-cells are the points of $A$, and $1$-cells are the arcs of $\gamma$
between these points. The preimage of this cell decomposition under $p$
is a cell decomposition of $X$. It essentially coincides with
the cell decomposition
described above, except that now we defined precisely which points
are vertices of multiplicity $2$, these are exactly the simple
preimages of the elements of the set $A$.
  
In classical literature
they usually consider the dual graph of this cell decomposition of
$X$, called the {\em line complex},
which can be obtained in the following way. Let $e_k$ be the edges
of $\gamma$ between adjacent pairs of the points of $A$.
Choose two point $x$ and $o$, one
in each complementary region of $\gamma$ and connect 
these points by simple disjoint arcs $\gamma_k$, so that $\gamma_k$
intersects $\gamma$ only once, and this intersection occurs at
a point of $e_k$. The points $x$ and $o$ and the arcs $\gamma_k$ form
a graph with two vertices and $q=\card A$ edges.
The preimage of this graph under $p$ is a properly embedded connected graph
in $X$ which is called the line complex of $(X,p)$. Every vertex of
the line complex has degree $q$ and every edge connects two vertices of
different type, one projecting to $x$ another to $o$. All possible line
complexes are characterized by these two properties.
The components of the complement of a line complex are called
faces.  There can be faces with an even number $2n$ of boundary edges,
they correspond to critical points of order $n$; and there can be 
unbounded faces with infinitely many boundary edges, they correspond
to the ``logarithmic branch points'' of the surface $(X,g)$.
The logarithmic branch points constitute the simplest
example of transcendental singularities.

Teichm\"uller's theorem mentioned in Section 1 implies that the conformal type
of $(X,g)$ is determined by its line complex.
 
The following question was asked
by A.~Epstein. Let $(\C,f)$ be a parabolic surface of Speiser class,
such that $f$ is of finite order. If $\psi$ is a homeomorphism
of the Riemann sphere, then $(\C,\psi\circ f)$ is parabolic because
as we mentioned before, Speiser surfaces have stable type.
Let $\phi$ be the uniformizing map of this deformed surface,
so that $g=\psi\circ f\circ\phi$ is a meromorphic function
in the plane. It is easy to show that $g$ also has finite order.
The question is whether the orders of $f$ and $g$ are the same.
In other words, does the line complex determine the order of
a function?
K\"unzi \cite{Kunzi} showed that this is not necessarily so for meromorphic
functions $f$, but the question remains unsolved for entire
functions\footnote{It is now solved in \cite{Bi}.}.
\vspace{.1in}

Next theorem gives
a complete characterization
of meromorphic functions whose inverses have finitely many singularities.
\vspace{.1in}

\noindent
4.2 {\sc Theorem}. {\em For a surface $(X,p)$ spread over the sphere,
the following conditions are equivalent:
\newline
(a) The set of transcendental singularities of $(X,p)$ is
finite, and $p$ has finitely many critical points. 
\newline
(b) $(X,p)$ is equivalent to $(\C,f)$ where $f$ is a meromorphic
solution of the
differential equation
\begin{equation}
\label{D}
{f'''}/{f'}-({3}/{2})\left({f''}/{f'}\right)^2=P,
\end{equation}
where $P$ is a rational function.

If $p$ has no critical points then $P$ is a polynomial,
and for every polynomial $P$ all solutions of $(\ref{D})$
are meromorphic functions in the plane.
The degree of $P$ is $n-2$, where $n$ is the number
of transcendental singularities of $(X,g)$.}
\vspace{.1in}

The case when $P$ is a polynomial and $p$ has no critical points
is due to R.~Nevanlinna \cite{Nev2}, and the generalization
with finitely many critical points to his student Elfving \cite{Elf}.
Every solution $f$
of the Schwarz differential
equation in (b) is a ratio of two linearly independent solutions
of the linear differential equation
\begin{equation}
\label{dd}
w''+(P/2)w=0.
\end{equation}
This provides very precise information on the asymptotic behavior
of $f$ sufficient to prove that (b) implies (a).

We give a sketch of Nevanlinna's argument, assuming for simplicity
that there are no critical points.
To show that (a) implies (b), Nevanlinna 
uses the line complex of $(X,p)$ to construct a sequence
of rational functions $(f_m)$ such that the surfaces
$(D_m,f_m),$ where $D_1\subset D_2\subset\ldots\to\C$
form an exhaustion
of the plane by some simply connected regions. These
surfaces $(D_m,f_m)$ 
approximate $(X,p)$ in the following sense.
There is an exhaustion $X_1\subset X_2\subset\ldots\to X$
by open topological
discs, such that each $X_m$ is isometric to a subset of
the surface $(D_m,f_m)$.
The rational functions $f_m$ have critical points
of high multiplicity,
one for each transcendental singularity of $(X,p)$.
So the number of
critical points of $f_m$ is bounded independently
of $m$.
The rational functions $f_m$
are normalized by $f_m(0)=0$ and
$f_m^\prime(0)=1.$
It follows from the Caratheodory Convergence theorem
that the sequence $(f_m)$ converges in some neighborhood
of $0$ to a holomorphic function $f$.
The crucial observation is that the Schwarzian derivatives
$$P_m=f_m^{\prime\prime\prime}/f_m^\prime-
(3/2)(f_m^{\prime\prime}/f_m^\prime)^2$$
are rational functions whose degrees
are bounded independently
of $m$. This is because the poles of
a Schwarzian derivative can occur only
at the critical points of the function,
and all these poles are of multiplicity two.
As the $P_m$ converge
in a neighborhood of $0$, they have to converge in the whole plane
to a rational function $P$. The limit function $f$ satisfies
the differential equation (\ref{D}), and has no critical points.
So $P$ is a polynomial, and $f$ has an analytic continuation
to a meromorphic function in the plane as a solution
of the differential equation (\ref{D}). Now it is easy to see that
$f$ is the meromorphic function associated with $(X,p)$.

To show that that (b) implies (a) we first notice that
$f$ has no critical points (each critical point in $\C$ is a pole
of the Schwarzian derivative, and we assume that $P$ is
a polynomial). Then one writes $f=w_1/w_0$
where $w_1$ are entire functions, linearly independent solutions
of the differential equation (\ref{dd}). Asymptotic integration
of this linear differential equation shows that the surface associated
to $f$ has finitely
many singularities, which proves (a). 

The singularities of $(\C,f)$ in Theorem 4.2 are simply related to the
so-called Stokes multipliers of the equation (\ref{dd}).
As it is easy to construct a surface as in (a) with prescribed
projections of singularities, one can derive the existence
of a linear differential equation (\ref{dd}) with prescribed
Stoke's multipliers \cite{Sibuya}.
\vspace{.1in}

In the case of infinitely many algebraic singularities, one cannot obtain
such complete results, but still for many subclasses of surfaces
of the Speiser class (which was defined in Section 1)
one can obtain
rather complete information about the asymptotic behavior of the uniformizing
functions \cite{DW,Gold,GO,Kunzi,Wit}.
\vspace{.1in}

{\bf 5.}
Our next example of an exact correspondence between a class of
functions and a class of surfaces involves a
class of cellular surfaces spread over the plane
and having a symmetry property. Suppose that an anticonformal involution
$s:X\to X$ is given. The set of fixed points of $s$ will be called the axis.
We say that $(X,p)$ is {\em symmetric} if $p\circ s=\overline{p}$ where
the bar stands for the complex conjugation. A symmetric surface spread over the
plane is called a {\em MacLane surface} if all its
critical points belong to the axis, and for
each transcendental singularity over $\C$,
the axis can serve
as an asymptotic curve.
Evidently, there can be
at most two transcendental singularities over $\C$.

When $X$ is identified with one of the standard domains, $\C$ or $\U$,
one can always ensure that the symmetry axis is $\R$ or $(0,1)$.
So holomorphic functions corresponding to MacLane surfaces are
real functions whose all singular values are real,
in other words, they are cellular functions with $\gamma=\R$.

It is easy to understand the nature of the cell decompositions
that correspond to MacLane functions. 


MacLane surfaces can be completely described by their singular
values. Suppose for simplicity that the sequence of critical points
on the axis is unbounded in both directions,
and that all critical points are simple.
Then the preimage of the
real axis under $p$ consists of the axis of symmetry and infinitely
many disjoint simple curves crossing the axis at the critical points
and tending to infinity in both directions. These curves together with
the axis form the $1$-skeleton of the cell decomposition of the plane.
The sequence
$(c_k)_{k\in\Z}$
of critical values, where each critical value is repeated according
to its multiplicity satisfies
$$(c_{k+1}-c_k)(c_k-c_{k-1})\geq 0,$$
and every sequence of critical values with this property can occur.

The corresponding class of functions is related to entire functions
of Laguerre--P\'olya (LP) class:
these are the real entire functions which are
limits of real polynomials with real zeros.
According to Laguerre and P\'olya,
this class LP has the parametric representation:
\begin{equation}
\label{lp}
z^m\exp(-a z^2+bz+c)\prod_k \left(1-\frac{z}{a_k}\right) e^{z/a_k},\quad
a\geq 0,\; b,c,a_k\in\R.
\end{equation}
It is easy to see that that LP-functions constitute a subclass of
MacLane's class.
In the case that the sequence critical values is infinite in both directions,
the critical values of an LP-function satisfy one additional restriction:
their signs alternate, 
\begin{equation}
\label{co}
c_{k+1}c_k\leq 0.
\end{equation}

\noindent
5.1 {\sc MacLane's Theorem}. {\em Every surface of MacLane's class
is parabolic. The derivatives of the corresponding entire functions 
constitute the class LP.}
\vspace{.1in}

This was proved in \cite{MacLane}. A very illuminating geometric proof
is given in \cite{Vinberg}. We sketch Vinberg's proof,
restricting ourselves for simplicity to the case when the sequence
of zeros is infinite in both directions. First consider MacLane surfaces
that satisfy the additional condition (\ref{co}).

Let $\Omega$ be the domain obtained from the plane by deleting
the vertical slits $L_k=\{ k\pi+i(y_k-t):0\leq t<\infty\}$, $k\in\Z$,
where $y_k$ are real numbers or $y_k=-\infty$ in which case the
slit is empty. We call such domains $\Omega$ {\em comb domains}.
We also require that the set $\{ k:y_k>-\infty\}$
be unbounded in both directions. Let $\theta:H\to\Omega$ be
a conformal map from the upper half-plane onto $\Omega$,
$\theta(\infty)=\infty$. Then the function
$f(z)=\exp(-i\theta(z))$, initially defined in 
the upper half-plane, can be extended by symmetry to an entire
function. It is easy to see that 
it belongs to MacLane's class,
with critical values $c_k=(-1)^ke^{y_k}$. 
By approximating
the region $\Omega$ by its 
intersections with vertical strips
$\{ z:|\Re z|<\pi m\}$, 
we obtain a sequence of real polynomials with real
zeros that converge to $f$. Thus $f$ is in LP class.

To prove that all functions of LP class satisfying
(\ref{co}) arise this way,
we just notice that every real polynomial of degree $2m$ with all zeros real
has a ``comb representation'' $p(z)=\exp(-i\theta(z))$
where $\theta$ is a conformal map of the upper half-plane
onto a comb domain in the strip $\{ z:|\Re z|<\pi m\}$. Passing to
the limit we obtain a comb representation of a given LP function satisfying (\ref{co}).

Thus we obtain a correspondence between functions
of the class LP and a subclass of MacLane surfaces characterized by
condition (\ref{co}). This correspondence becomes bijective if
we factor sequences of critical values by shifts, and
entire functions by the change of the independent variable 
$az+b,$ $a>0, b\in\R$.

One can even tell
explicitly in terms of $c_k$ when $a=0$ in (\ref{lp}),
see \cite{GoldbergCV}, \cite{ASWE}.

To prove the general case of MacLane's theorem,
one can do an approximation argument similar to that in
the proof of Theorem 4.2. The crucial fact here is
the following. 
\vspace{.1in}

{\em Let $(p_n)$ be a sequence of real polynomials whose all
zeros are real, and assume that $(p_n)$ converges in a
neighborhood of $0$ to a non-zero function. Then the sequence
converges on every compact subset of the plane.}
\vspace{.1in}

Let $(X,p)$ be a surface of MacLane class, and $f=p\circ\phi$
the corresponding function (which is \`a priori holomorphic
either in $\C$ or in $\U$). Using the cell decomposition,
one can find an exhaustion of $X$ by open topological discs
$X_1\subset X_2\subset\ldots\to X$ such that each $X_k$ is isometric
to a subset of a surface corresponding to a real polynomial $f_k$
with all critical points real. By Caratheodory's theorem,
these polynomials converge to $f$ in a neighborhood of $0$.
 Then the derivatives $f_k^\prime$ are polynomials
with all zeros real, and they converge in a neighborhood of $0$
to $f'$. Applying the above proposition we conclude that $f_k^\prime$
converge in the whole plane, and thus $f_k$ also converge in
the whole plane. So $f$ is entire.
\vspace{.1in}

The class LP and its geometric
characterization occur in many 
questions of analysis. 

The most striking application of the geometric characterization
of the class LP is to the spectral theory of self-adjoint
second order
differential operators on the real line with periodic potentials (Hill
operators).
It was discovered by Krein \cite{Krein} that Lyapunov functions of
periodic strings are exactly those real entire functions $f$
of genus zero, with positive roots, which have the property that
all solutions $z$ of the equation $f^2(z)=1$ are non-negative.

These functions constitute a subclass
of LP which can be explicitly described
in terms of their comb domains. 

We include a brief explanation of this connection.
Consider a second-order operator
$$L(y)=-y^{\prime\prime}+v(x)y,$$
where $v$ is a periodic function of period $\pi$.
Then the spectra of the boundary value problems
with periodic or antiperiodic boundary conditions
are known to be real and bounded from below.

Let $c(x,\lambda)$ and $s(x,\lambda)$ be solutions of
the differential equation $L(y)=\lambda y$ with
the following initial conditions:
$$s(0,\lambda)=c'(0,\lambda)=0\quad s'(0,\lambda)=c(0,\lambda)=1.$$
One can show (by reducing the differential equation to an integral
equation which can be solved by the standard iteration procedure)
that $c$ and $s$ are entire functions of $\lambda$ of order at most $1/2$,
normal type. Following Lyapunov, we are trying to find out
when all solutions of our differential equation are
bounded on the real line. For this purpose, we consider the monodromy
matrix
$$T(\lambda)=\left(\begin{array}{ll}c(\pi,\lambda)&s(\pi,\lambda)\\
c'(\pi,\lambda)&s'(\pi,\lambda)\end{array}\right),$$
whose determinant is $1$ because of the Liouville's theorem and the
initial conditions. One half of the trace of $T$ is called
the {\em Lyapunov function},
$$u(\lambda)=(c(\pi,\lambda)+s'(\pi,\lambda))/2.$$
This is a real entire function of order at most $1/2$.
The eigenvalues $z_1$ and $z_2$
of the monodromy matrix are the solutions of the quadratic equation
$z^2-2u(\lambda)z+1=0$.
It is easy to see that all solutions of the differential
equation will be bounded if $u^2(\lambda)\leq 1$.
Thus the real line is divided into two sets:
the stability zone where $u^2\leq 1$ and the instability zone.
On the boundary of the stability zone, we have $u(\lambda)=\pm 1$.
Let us consider this last equation in the complex plane. If $u(\lambda)=\pm 1$,
then
$T(\lambda)$ has one multiple eigenvalue $\pm 1$. This means
that our differential equation $L(y)=\lambda y$ has either periodic
solution of period $\pi$ or an ``antiperiodic'' solution
$y(0)=-y(\pi)$. As both periodic and antiperiodic boundary value problems
are self-adjoint, all their eigenvalues should be real,
and we conclude that the equation $u^2(\lambda)=1$ has only real
solutions.

So both $u-1$ and $u+1$ belong to the Laguerre--P\'olya class
(as real entire functions of order less than two, with real zeros).
So all zeros of $u'$ are also real, so $u$ belongs to
the MacLane class. Now it is easy to conclude that in fact $u\in LP$.

For Lyapunov functions, Marchenko and Ostrovski derived a different
version of comb representation, namely $u(z)=\cos\theta(z)$,
where $\theta$ is a conformal map of the upper half-plane
onto a region obtained from the upper half-plane by deleting vertical
slits from the points $\pi k$ of height $h_k$. Thus we have
\vspace{.1in}

\noindent
5.2 {\sc Theorem of Marchenko and Ostrovskii} {\em Let $f$ be a real
entire function such that the equation $f^2(z)=1$ has only real
solutions. Then $f=\cos\theta,$ where $\theta$ is the conformal map
we just described.}
\vspace{.1in}

The theorem does not contain any a priori growth assumptions,
so one has to show first that $f$ is of order at most $2$.
This can be done by elementary Nevanlinna theory, but 
for applications to Lyapunov functions such generalization is unnecessary.

Krein proved the converse theorem, that every real entire function
of genus zero such that all solutions of $f^2(z)-1=0$ are real,
is a Lyapunov function of a periodic string.
The ``string'' in this result can be quite singular, the mass of
the string is in general not a function but an arbitrary measure.
So the question of characterizing Lyapunov functions arises for
various subclasses when one impose conditions on the potential $v$,
for example $v\in L^2$.

Developing Krein's ideas, Marchenko and Ostrovskii \cite{MO}
obtained a
pa\-ra\-met\-ri\-za\-tion of self-adjoint
periodic Hill operators with $L^2$ potentials
in terms of their spectral data.
In a recent paper \cite{Tkachenko}, Tkachenko extended this result to
some non self-adjoint Hill operators (with complex potentials) by
considering small perturbations of MacLane surfaces which are no longer
symmetric, and establishing an exact correspondence between a class of entire
functions and a class of surfaces in the spirit of Theorem~5.1.

We finish this section by mentioning a classical problem about real entire
functions with real zeros which was recently solved. 

It is evident from the definition that all derivatives of a function
of the class LP have only real zeros.
The converse is also true, and in a very
strong sense:
if $f$ is a real entire function, and all zeros of $ff''$ are real, then
$f\in$ LP \cite{BEL}. This was conjectured by Wiman in 1911, and the final
result completes a long line of development with
important contributions of Levin, Ostrovskii, Hellerstein, Williamson
and Sheil-Small. Recently Edwards and Hellerstein (for the case of finite
order) and
Langley (for the case of infinite order)
extended the result by replacing $f''$ by $f^{(k)}$
with any $k\geq 2$.

\vspace{.1in}

{\bf 6.} 
Here we mention a partial
generalization of Theorem 5.1 to meromorphic functions.
Consider the class of meromorphic functions of the form
\begin{equation}
\label{3}
f(z)=e^{\sigma z}\frac{\prod(1+z/a_k)}{\prod(1-z/b_k)}.
\end{equation}
where $\sigma\geq 0$ and $(a_k)$ and $(b_k)$ are two
increasing sequences of positive numbers, finite or infinite
(possibly empty), and such that
$$\sum\frac{1}{a_k}+\sum\frac{1}{b_k}<\infty.$$
Let us denote by $(x_k)$ the sequence of real critical points of $f$.
They are enumerated by positive integers on the positive ray and
by negative integers on the negative ray.
We assume that this sequence is increasing, each point is
repeated according to its multiplicity, and $x_{-1}<0<x_1.$
Let $c_k=f(x_k)\in\bR=\R\cup\{\infty\}$
be the corresponding sequence of critical values.
\vspace{.1in}

\noindent
6.1 {\sc Theorem}. {\em All critical points of a function $f$
of the form $(\ref{3})$ are real. For a sequence $(c_k)$
to be a sequence of critical values of the function of
the form $(\ref{3})$, it is necessary and sufficient that
the following conditions are satisfied:
\newline
$(-1)^kc_k\geq 0$,
\newline
if $c_k=0$ then $0\in\{ c_{k-1},c_{k+1}\}$; if $c_k=\infty$ then
$\infty\in\{c_{k-1},c_{k+1}$, and
\newline
$|c_{-k}|<|c_k|$ for all $k>0$ for which both $c_k$ and $c_{-k}$ are
defined.}
\vspace{.1in}

The class of meromorphic functions of the form (\ref{3})
has another interesting parametrization: Taylor coefficients
of these functions are exactly the totally positive sequences
\cite{Aissen}.

The class of functions (\ref{3}) does not exhaust all real
meromorphic functions
whose critical points are real, even if we restrict ourselves
to functions without asymptotic values.

Thus no satisfactory analog of Theorem 5.1 is known
for meromorphic functions and open surfaces spread over the sphere.
However, there is a related result for {\em rational functions} which was
used in \cite{EG} to prove a special case of an intriguing conjecture
in real algebraic geometry.
\vspace{.1in}

\noindent
6.2 {\sc Theorem}. {\em If all critical points of a rational function $f$
belong
to a circle $C$ on the Riemann sphere, then $f(C)$ is a subset of
a circle.}
\vspace{.1in}

Let us call two rational functions $f$ and $g$ equivalent if $f=\ell\circ g$
where $\ell$ is a fractional linear transformation. Equivalent 
functions have the same critical points. We may assume
without loss of generality that the circle $C$ in Theorem 6.2 is the real
line. Then Theorem 6.2 says that whenever the critical points of a rational
function are real, it is equivalent to a real rational function.
A rational function with prescribed critical points can be obtained as
a solution of a system of algebraic equations, so Theorem 6.2 implies that
all solutions of this system are real whenever the coefficients are real.
We will see in a moment the geometric significance of this system of
algebraic
equations.

It is very easy
to prove Theorem 6.2 for polynomials,
because every two polynomials with
the same critical points are equivalent (which means that solution of
the system
of algebraic equations mentioned above is essentially unique in this case).
This is not so for rational
functions: it turns out that for given $2d-2$ points in general position
on the sphere
there are finitely many, namely
$$u_d=\frac{1}{d}\binom{2d-2}{d-1},\quad\mbox{the $d$-th Catalan number,}$$
of classes of
rational functions of degree $d$ which share these critical points. This is due
to L. Goldberg \cite{LGold} who reduced the problem to the following classical
problem of enumerative geometry: given $2d-2$ lines in general
position in (complex) projective space,
how many subspaces of codimension $2$ intersect
all these lines? The answer to this last problem, the Catalan number
$u_d$, was obtained by Schubert
in 1886 who invented what is now known as ``Schubert Calculus'' to solve this
and similar enumerative problems. 

Notice that rational functions exhibit a very special property
in Goldberg's result: as we mentioned above, there is only one class of
polynomials
with prescribed critical points, and similarly there is only
one class of Blaschke products with prescribed critical points \cite{Zak}.

The general question of how many solutions to a problem of enumerative geometry
can be real was asked by Fulton in \cite{Fulton}. A specific conjecture
about the problem of finding subspaces of appropriate
codimension intersecting given real subspaces was made by
B. and M. Shapiro. For the problem of Schubert calculus stated above,
this conjecture says that
if all $2d-2$ given lines are tangent to the rational normal curve
$z\mapsto(1:z:\ldots:z^d)$ at real points,
then all $u_d$ subspaces of codimension $2$ which intersect
these lines are real (can be defined by real equations). This statement is
equivalent to Theorem 6.2.

The proof of Theorem 6.2 is based on an explicit description
of surfaces spread over the sphere which correspond to real rational functions
with real critical points in the spirit of Vinberg's work about
the MacLane's class mentioned in section 4.

Let $R$ be the class of {\em real} rational functions $f$
whose all critical
points are real and simple.
To describe the Riemann surface of $f^{-1}$ 
we consider the {\em net} $\gamma_f$ which is
the preimage $\gamma_f=f^{-1}(\R\cup\infty)$ 
modulo homeomorphisms of the Riemann sphere commuting with complex
conjugation. A net consists of simple arcs which meet only at the critical
points of $f$. To each of these arcs we prescribe a {\em label} equal to
the length of its image under $f$. One can describe explicitly all
labeled nets which may occur from this construction.
It turns out that labeled nets give a parametrization of the class $R$ of
rational functions. This parametrization has an advantage that
it clearly separates the discrete, topological parameter (the net)
from the continuous parameters (the labeling), and it turns out that
the set of possible labelings
of a given net has simple topological structure:
it is a convex polytope.

Using topological methods, we show in \cite{EG} that
for every net and for every set of $2d-2$ points on the real line
there exists a rational function of the class $R$ with this net
and these critical points. On the other hand, a simple combinatorial
argument shows that the number of possible nets on $2d-2$ vertices is
equal to the Catalan number $u_d$. Thus one obtains $u_d$ classes of
real rational functions of degree $d$ with prescribed real critical points.
Comparison with L. Goldberg's
result shows that in fact we constructed {\em all}
classes of rational functions with prescribed real critical points. Thus
all such classes contain real functions.

Since the publication of \cite{EG}, two new proof of Theorem 6.1
appeared. One is a very much simplified version of
the original proof; it still uses the parametrization of the
class $R$ described above, but avoids the Uniformization theorem \cite{EG2}.
Another proof \cite{MTV} is based on completely different ideas,
coming from the theory
of exactly solvable models of statistical mechanics. It proves
not only 
Theorem 6.1 but also its multi-dimensional generalization
to rational curves in projective spaces.

The author thanks Mario Bonk, David Drasin, Juha Heinonen, Pietro Poggi-Coradini
and Misha Sodin for help and encouragement.

\end{document}